\documentclass[10pt]{amsart}
\usepackage{setspace,amscd}%

\author[D. Khosla]{Deepak Khosla}
\address{Department of Mathematics \\  University of Texas at Austin \\
  1 University   Station C1200 \\ Austin, Texas 78712 \\ USA}
\email{dkhosla@math.utexas.edu}

\date{\today}

\title[Moduli of Curves with Linear Series]
{Moduli Spaces of Curves with Linear Series and the Slope Conjecture}

\newtheorem{thm}{Theorem}[section]
\newtheorem{lemma}[thm]{Lemma}
\newtheorem{cor}[thm]{Corollary}
\newtheorem{prop}[thm]{Proposition}

\theoremstyle{definition}
\newtheorem{defn}[thm]{Definition}
\newtheorem{rmk}[thm]{Remark}
\newtheorem{notation}{Notation}

\newcommand{\oh}[2]{\ensuremath{\mathcal{O}_{#1}(#2)}}

\DeclareMathOperator{\Sym}{Sym}
\newcommand{\tensor}{\otimes}
\newcommand{\caniso}{\simeq}
\newcommand{\isom}{\cong}
\newcommand{\ds}{\oplus}
\newcommand{\Ds}{\bigoplus}
\newcommand{\dual}[1]{{#1}^{\vee}}

\newcommand{\coh}[2]{H^{#1}(#2)}

\newcommand{\incl}{\hookrightarrow}
\newcommand{\union}{\cup}
\newcommand{\Union}{\bigcup}

\newcommand{\Z}{\mathbf{Z}}
\newcommand{\Q}{\mathbf{Q}}
\newcommand{\C}{\mathbf{C}}
\newcommand{\Proj}{\mathbf{P}}
\newcommand{\Grass}{\mathbf{G}}
\newcommand{\hilb}[1]{\Hilb_{\Z}^{#1}\Proj^2}
\newcommand{\sch}[1]{\sigma_{#1}}
\newcommand{\ls}[1]{\lvert{#1}\rvert}
\newcommand{\st}[2]{\left\{#1\mid #2\right\}}
\DeclareMathOperator{\Pic}{Pic}
\DeclareMathOperator{\Bl}{Bl}
\DeclareMathOperator{\Hilb}{Hilb}

\DeclareMathOperator{\ch}{ch}
\DeclareMathOperator{\td}{td}

\newcommand{\sC}{\mathcal{C}}
\newcommand{\sL}{\mathcal{L}}
\newcommand{\sV}{\mathcal{V}}

\newcommand{\g}[2]{\mathfrak{g}_{#2}^{#1}}
\newcommand{\grd}{\g{r}{d}}
\newcommand{\Grd}{G^r_d}

\newcommand{\sgrd}{\mathcal{G}^r_d}
\newcommand{\scrd}{\mathcal{C}^r_d}
\newcommand{\gsix}{\mathcal{G}^6_{24}}

\newcommand{\csix}{\mathcal{C}^6_{24}}
\newcommand{\gsixm}{\gsix}

\newcommand{\Mgb}{\overline{\mathcal{M}}_g}
\newcommand{\Mgnb}{\overline{\mathcal{M}}_{g,n}}
\newcommand{\Mg}{\mathcal{M}_g}
\newcommand{\Mtl}{\widetilde{\mathcal{M}}_{g,1}}
\newcommand{\Ctl}{\widetilde{\mathcal{C}}_{g,1}}
\newcommand{\Mttl}{\widetilde{\mathcal{M}}_{2,1}}
\newcommand{\Mogb}{\overline{\mathcal{M}}_{0,g}}
\newcommand{\mtl}{\widetilde{\mathcal{M}}}
\newcommand{\m}{\mathcal{M}}
\newcommand{\mb}{\overline{\mathcal{M}}}
\newcommand{\mirr}{\m_{21}^{\mathrm{irr}}}
\newcommand{\cb}{\overline{\mathcal{C}}}

\newcommand{\ac}{\alpha}
\newcommand{\bc}{\beta}
\newcommand{\cc}{\gamma}
\newcommand{\sg}{\sigma}
\newcommand{\Sg}{\Sigma}
\newcommand{\Sgt}{\overline{\Sigma}}
\newcommand{\lam}{\lambda}
\newcommand{\de}{\delta}
\newcommand{\De}{\Delta}
\newcommand{\dn}{\delta_0}
\newcommand{\di}{\delta_i}
\newcommand{\eps}{\epsilon}
\newcommand{\epsi}{\epsilon_i}
\newcommand{\tE}{\widetilde{E}}
\newcommand{\tD}{\widetilde{D}}

\newcommand{\mess}{\xi}
\newcommand{\chow}[2]{A^{#1}({#2})}
\newcommand{\pf}{\eta_*}
\newcommand{\pls}{\pi_*}

\begin{document}

\bibliographystyle{hamsplain}

\begin{abstract}
  We describe the moduli space $\sgrd$ of triples consisting of a
  curve $C$, a line bundle $L$ on $C$ of degree $d$, and a linear
  system $V$ on $L$ of dimension $r$. This moduli space extends over a
  partial compactification $\mtl_g$ of $\m_g$ inside $\Mgb$. For the
  proper map $\eta\colon\sgrd\to \mtl_g$, we compute the push-forward
  on Chow 1-cocyles in the case where $\eta$ has relative dimension
  zero. As a consequence we obtain another counterexample to the
  Harris-Morrison slope conjecture as well as an infinite sequence of
  potential counterexamples.
\end{abstract}

\maketitle{}

\section{Introduction}
\label{sec:introduction}

One of the fundamental objects in the study of the birational geometry
of a projective variety is its cone of effective divisors. 
Deligne and Mumford \cite{DelMum} constructed a compactification
$\Mgb$ of the moduli space $\Mg$ of smooth genus-$g$ algebraic curves. By a
theorem of Harer \cite{Harer} in topology, for $g\ge 3$,
\begin{equation*}
  \Pic \Mgb\tensor\Q =
  \Q\lam\ds\Q\dn\ds\Q\de_1\ds\dotsb\ds\Q\de_{\lfloor g/2 \rfloor},
\end{equation*}
where $\lam$ comes from a class on $\Mg$ and the $\di$ are boundary
classes. It is natural to ask the question: which linear combinations
of $\lam$ and the $\di$ are effective? If we consider only classes of
the form $a\lam - b\de$, where
\begin{equation*}
  \de = \dn + \de_1 + \dotsb + \de_{\lfloor g/2 \rfloor}
\end{equation*}
is the total boundary class, then the shape of the effective cone on
this plane is determined by the single number
\begin{equation*}
  s_g = \inf 
    \st{\frac{a}{b}}{a,b\ge0\text{, } a\lam - b\de\text{ is effective}},
\end{equation*}
which is known as the \emph{slope} of $\Mgb$. By explicitly constructing
effective divisors with small slope, Harris and Mumford \cite{HarMum}
deduced that $\Mg$ is of general type when $g\ge 24$. Conversely, a
suitable lower bound on the slope of $\Mgb$ would force the Kodaira
dimension of $\Mg$ to be $-\infty$. In 1990 Harris and Morrison
\cite{HarMor} conjectured that
\begin{equation*}
  s_g \ge 6 + \frac{12}{g+1}.
\end{equation*}

Recently Farkas and Popa \cite{FarkasPopa} found a counterexample to
this conjecture in genus 10. Specifically they considered the divisor
\begin{equation*}
  D = \st{[C]\in\m_{10}}{C\text{ lies on a K3 surface}}
\end{equation*}
and showed that its class is
\begin{equation*}
  7\lam - \dn - 5\de_1 - \dotsb .
\end{equation*}
This has a slope of 7, which is less than $6+\frac{12}{11}$. Farkas
has gone further and shown that $\Mg$ is of general type when $g$ is
22 or 23 \cite{Farkas.Syzygies}.
Though the divisor $D$ appears to be an isolated example, in fact, it may be
re-written as the divisorial component of
\begin{quote}
  the locus of $[C]\in\m_{10}$ for which there is an embedding $C\incl\Proj^4$
      of degree 12 such that $C$ lies on a quadric.
\end{quote}
This suggests a host of generalizations. For example, we may define
$D^{r,k}_{d,g}$ to be
\begin{quote}
   the locus of $[C]\in\m_{g}$ for which there is an embedding $C\incl\Proj^r$
      of degree $d$ such that $C$ lies on a hypersurface of degree $k$.
\end{quote}
The problem is then to determine when such cycles are divisors and, in
those cases, to compute their classes. Our approach is to consider the
moduli space $\sgrd(\Mg)$ of \emph{curves with linear series}. We set
\begin{equation*}
  \sgrd(\Mg)=\sgrd=\st{(C,L,V)}{[C]\in\Mg,L\in\Pic^d C,V\subset\coh{0}{L}},
\end{equation*}
where $V$ is a subspace of dimension $r+1$. Then $D_{d,g}^{r,k}$ is
naturally the image of the subscheme $\tD$ in $\sgrd$,
\begin{equation*}
  \tD = \st{(C,L,V)}{\Sym^k V\to\coh{0}{L^{\tensor k}}\text{ has a kernel}}
\end{equation*}
If it is of the expected dimension, the class of $\tD$ is then easily
evaluated in terms of certain standard classes on $\sgrd$. It remains
to compute the proper push-forward map for the morphism
\begin{equation*}
  \eta\colon\sgrd\to\Mg .
\end{equation*}
In our main result, Theorem \ref{main-thm}, we compute $\eta_*$ on
divisors in the case where $\eta$ has relative dimension 0.

In Section \ref{sec:gener-brill-noeth} we work out the case where
$g=21$, $r=6$, $d=24$, and $k=2$, so that the cover
\begin{equation*}
  \gsix \to \m_{21}
\end{equation*}
is generically finite, and the condition of lying on a quadric imposes
one condition on a $\g{6}{24}$. We show that $D_{24,21}^{6,2}$ has a
divisorial component and use Theorem \ref{main-thm} to compute its
class. It turns out that $D_{24,21}^{6,2}$ provides another
counterexample to the Slope Conjecture. In fact, one may consider
\begin{equation*}
  (g,r,d) =\bigl( m(2m+1),2m,2m(m+1)\bigr)
\end{equation*}
for any integer $m\ge1$. In this case the Brill-Noether number $\rho$
is zero, and it is one condition to lie on a quadric. Although at this
time we do not know how to prove that the $D_{d,g}^{r,k}$ are divisors
for all $m$, our calculations show that if they were, they would all
provide counter-examples to the Slope Conjecture.

In Section~\ref{sec:moduli-curves-with} we outline the definition and
basic properties of the moduli space $\sgrd$ before stating the main
theorem. In Section~\ref{sec:spec-famil-curv} we give the statements
of a series of calculuations over special families of stable curves.
These calculations assemble to give the main result. Finally
Section~\ref{sec:proofs-lemmas} is devoted to the proofs of the lemmas
stated in Section~\ref{sec:spec-famil-curv}.

This work was carried out for my doctoral thesis under the supervision
of Joe Harris. I would like to thank Ethan Cotterill, Gavril Farkas,
Johan de Jong, Martin Olsson, Brian Osserman, and Jason Starr for
helpful conversations.

\begin{notation}
  Unless otherwise indicated, all schemes are of finite type over
  $\C$. For a scheme (or Deligne-Mumford stack) $X$ of dimension $n$, we
  write $A^k(X)$ for the Chow group $A_{n-k}(X)$.
\end{notation}

\section{A Generalized Brill-Noether Divisor}
\label{sec:gener-brill-noeth}

In this section we work on $\mirr$, the locus of irreducible stable
genus-21 curves, and we let $\gsix$ denote $\gsix(\mirr)$. 
This space may be informally defined as the locus of triples
\begin{equation*}
  \st{(C,L,V)}{[C]\in\mirr,L\in\overline{\Pic}^{24} C,V\subset\coh{0}{L}},
\end{equation*}
where $L$ is allowed to be a torsion-free sheaf of rank 1.  In
Section~\ref{sec:moduli-curves-with} we will give a precise definition
of this space, which differs slightly from the above.

We first establish the following result.

\begin{lemma}
  \label{lem:irred}
  The space $\gsix(\m_{21})$ is irreducible.
\end{lemma}

\begin{proof}
  Note that a $\g{6}{24}$ is residual to a $\g{2}{16}$; that is
  \begin{equation*}
    L\in W^6_{24}(C) \iff K_C\tensor L^* \in W^2_{16}(C)
  \end{equation*}
  for any smooth curve $C$ of genus 21.
  Thus there is a dominant rational map
  \begin{equation*}
    V_{16,21} \longrightarrow \gsix
  \end{equation*}
  from the Severi variety of irreducible plane curves of degree 16 and
  genus 21. Since $V_{16,21}$ is irreducible \cite{SeveriProb} and
  maps dominantly to $\gsix$, so $\gsix$ is irreducible.
\end{proof}

If $\pi\colon\csix\to\gsix$ is the universal curve, we let
$\sL\to\csix$ be a universal line bundle and $\sV\subset\pi_*\sL$ the
universal rank-7 subbundle.

\begin{defn}
  Consider the open set $U\subset\gsix$ over which
  $\sL$ is a line bundle. Over $U$ there is a map $m$ of
  vector bundles of rank 28,
  \begin{equation*}
    m\colon\Sym^2\sV \to \pls\sL^{\tensor 2}
  \end{equation*}
  Let $\tE$ be the closure in $\gsixm$ of the singular locus of
  $m$. Since the complement of $U$ has codimension 2, the class of the
  degeneracy locus in $U$ extends uniquely to $\gsixm$.
\end{defn}

\begin{prop}
  \label{prop:E-tilde-divisor}
  The scheme $\tE$ has codimension 1 inside the irreducible component
  of $\gsix$ which dominates $\m_{21}$.
\end{prop}

\begin{proof}
  Since $\gsix(\m_{21})$ is irreducible, it suffices to exhibit a
  smooth curve with a $\g{6}{24}$ not lying on a quadric. Let
  \begin{equation*}
    S = \Bl_{21}\Proj^2
  \end{equation*}
  be the blow-up of $\Proj^2$ at 21 general points, and consider the
  linear system
  \begin{equation*}
    \nu = \Bigl\lvert 13H - 2\sum_{j=1}^9E_j - 3\sum_{k=10}^{21}E_k \Bigr\rvert
  \end{equation*}
  on $S$, where $H$ is the hyperplane class and $E_i$ are the
  exceptional divisors.
  A calculation in Macaulay (see the appendix) shows that
  a general member $C$ of $\nu$ is irreducible
  and smooth of genus 21. The series
  \begin{equation*}
    \Bigl\lvert 6H - \sum_{i=1}^{21}E_i\Bigr\rvert
  \end{equation*}
  embeds $S$ in $\Proj^6$ as the rank-2 locus of general $3\times6$
  matrix of linear forms \cite[Section 20.4]{Eisenbud.comm.alg}.
  The ideal of $S$ is therefore generated by cubics, so $S$ does not
  lie on quadric. It follows that $C$, which
  embeds in $\Proj^6$ in degree 24,
  does not lie on a quadric.

\end{proof}

\begin{defn}
  \label{def:E}
  Let $E$ be the effective codimension-1 Chow cycle which is the image
  of $\tE$ under the map
  \begin{equation*}
    \eta\colon \gsixm \to \mirr
  \end{equation*}
\end{defn}

To compute the class of $E$, we begin by expressing the class of $\tE$
in terms of simpler classes on $\gsix$. Let
\begin{align*}
  \ac &= \pls c_1(\sL)^2 \\
  \bc &= \pls c_1(\sL)\cdot c_1(\omega) \\
  \cc &= c_1(\sV)
\end{align*}

\begin{prop}
  \label{prop:class-E-tilde}
  The class of $\tE\subset\gsixm$ is
  \begin{equation*}
    2\ac - \bc + \lam -8\cc
  \end{equation*}
\end{prop}

\begin{proof}
  By Porteous, the class of $\tE$ is
  \begin{equation*}
    c_1(\pls\sL^{\tensor 2}) - c_1(\Sym^2\sV)
  \end{equation*}
  By Grothendieck-Riemann-Roch applied to the projection
  $\pi\colon\csix\to\gsix$ from the universal curve,
  \begin{equation*}
    \begin{split}
      \ch (\pls\sL^{\tensor2})&=
      \pls\bigl[\ch(\sL^{\tensor2})\cdot\td_{\csix/\gsix}\bigr] \\
      &= \pls\biggl[\Bigl(1+2c_1(\sL) + 2c_1(\sL)^2+\dotsb\Bigr)\biggr. \\
      & \qquad \qquad \left.\cdot \Bigl(1-\frac{c_1(\omega)}{2} + 
        \frac{c_1(\omega)^2+\kappa}{12}+\dotsb\Bigr)\right] \\
      &= 28 + (2\ac - \bc + \lam)+\dotsb
    \end{split}
  \end{equation*}
  where $\omega$ is the relative dualizing sheaf for $\pi$, and
  $\kappa$ is the divisor of nodes on $\csix$.
  Also,
  \begin{equation*}
    c_1(\Sym^2\sV) = 8c_1(\sV) = 8\cc\text{,}
  \end{equation*}
  and the proposition follows.
\end{proof}

\begin{prop}
  \label{prop:class-E}
  The class of $E\subset \mirr$ is given as
  \begin{equation*}
    [E] = 2459\lam - 377\dn \text{.}
  \end{equation*}
\end{prop}

\begin{proof}
  By Proposition
  \ref{prop:class-E-tilde} and Theorem \ref{main-thm},
  \begin{equation*}
    \begin{split}
      [E] = \pf[\tE] &= \pf(2\ac - \bc + \lam  - 8\cc) \\
      &= \frac{2459N}{95}\lam - \frac{377N}{95}\dn \text{,}
    \end{split}
  \end{equation*}
  where $N$ is the degree of $\eta$.
\end{proof}

\begin{cor}
  The slope conjecture is false in genus 21.
\end{cor}

\begin{proof}
  Since
  \begin{equation*}
    \frac{2459}{377} < 6 + \frac{12}{22},
  \end{equation*}
  this is an immediate consequence of \cite[Corollary 1.2]{FarkasPopa}.
\end{proof}

\begin{rmk}
  For any integer $m\ge1$, if we let
  \begin{equation*}
    (g,r,d) =\bigl( m(2m+1),2m,2m(m+1)\bigr) \text{,}
  \end{equation*}
  then, as mentioned in the introduction, we can consider an analogous locus
  \begin{equation*}
    \tE \subset \sgrd(\Mg^{\mathrm{irr}}) \text{.}
  \end{equation*}
  At this time, we do not have a result similar to
  Lemma~\ref{lem:irred} which allows us to show that $\tE$ is actually
  a divisor. Were this the case, however, identical computations to
  those above would show that the ratio $a/b_0$ of the coefficients of
  $\lam$ and $\dn$  in $E$ is less than $6 + 12/(g+1)$. Specifically,
  we compute the difference
  \begin{equation*}
    6 - \frac{12}{g+1} - \frac{a}{b_0} = 
    \frac{36{m}^{5}-24 {m}^{4}-57 {m}^{3}+48 {m}^{2}+3 m-6}
    {16 {m}^{9}-8 {m}^{8}-4 {m}^{7}-10 {m}^{6}+23 {m}^{4}+16
      {m}^{3}+13 {m}^{2}+2 m}  \text{.}
  \end{equation*}
\end{rmk}

\section{Moduli of Curves with Linear Series}
\label{sec:moduli-curves-with}

In this section we give the basic definitions needed to state the main
theorem. We begin by constructing a partial compactification of $\Mg$
inside $\Mgb$
over which we can extend the space $\sgrd(\Mg)$ of curves with linear series.

\begin{defn}
  \label{def:tree-like}
  For $i\in\{0,1,\dotsc,g\}$, let $B_i \subset \Mgnb$ be the locus of
  stable pointed curves which are the union of a smooth curve of genus
  $i$ and a smooth curve of genus $g-i-1$ meeting nodally at two
  points.  Define $\mtl_{g,n}$ to be the open substack of $\Mgnb$
  which is the complement of the closure of $\Union B_i$.
\end{defn}

\begin{defn}
  \label{def:script-Grd}
  We define the Deligne-Mumford stack $\sgrd(\mtl_{g,n})$ of curves with
  linear series as follows. Fiberwise, over each irreducible curve
  $[C]\in\mtl_{g,n}$, we consider torsion-free rank-1 coherent sheaves
  $L$ on $C$ together with a vector subspace of the space of global
  sections of $L$. Over reducible fibers, we consider limit linear
  series as defined by Eisenbud-Harris in \cite{LimLinSer}. These
  constructions have been made functorial by Altman-Kleiman
  \cite{AltmanKleiman1} and Osserman \cite{Osserman}. Osserman treats
  the case where there are at most two irreducible components and
  there is no monodromy among the components. The details of the
  general contruction will appear in a forthcoming paper
  \cite{limit.linear.series}. 
  An important fact is that given a reducible family of curves, the
  scheme structure on the open locus in $\sgrd$ of refined series
  coincides with the natural subscheme structure inside the product of
  Grassmannians.

  The resulting Deligne-Mumford stack ${\sgrd}'$ may not be
  irreducible, even over $\m_{g,n}$; we let $\sgrd=\sgrd(\mtl_{g,n})$
  be the unique irreducible component of ${\sgrd}'$ which dominates
  $\m_{g,n}$.

  The morphism
  \begin{equation*}
    \eta\colon \sgrd \to\mtl_{g,n}
  \end{equation*}
  is representable and proper. In the case where
  \begin{equation*}
    \rho = g - (r+1)(g-d+r)
  \end{equation*}
  is non-negative, $\eta$ is
  generically smooth of relative dimension $\rho$.
  Given $n$ ramification conditions $\beta_1,\dotsc,\beta_n$, the substack
  \begin{equation*}
    Z_{\beta_1,\dotsc,\beta_n} \subset \sgrd
  \end{equation*}
  of linear series with specified ramification at the marked points
  has the expected generic relative dimension
  \begin{equation*}
    \rho - \sum_i |\beta_i|
  \end{equation*}
  over $\mtl_{g,n}$. We say that a stable marked curve
  $(C,p_1,\dotsc,p_n)$ is \emph{Brill-Noether general} if, for all
  possible sets of ramification conditions $\beta_1,\dotsc,\beta_n$,
  the marked curve $(C,p_1,\dotsc,p_n)$ lies in the dense open over which
  $\eta|_{Z_{\beta_1,\dotsc,\beta_n}}$ is flat.
  In the following we work over $\Mtl$ in order to be able to
  consistently define the universal line and vector bundles.
  If $\pi\colon\scrd \to \sgrd$ is the universal curve, and
  $\sg\colon\sgrd\to\scrd$ is the marked section, then there is a
  universal coherent sheaf $\sL$ on $\scrd$, flat over $\sgrd$, with the
  properties that
  \begin{itemize}
  \item $\sL$ has torsion-free rank-1 fibers
  \item $\sL$ has degree $d$ on the component of each fiber which
    contains the marked point and has degree $0$ on all other components
  \item $\sL$ is trivialized along the marked section:
    $\sg^*\sL\caniso\mathcal{O}_{\sgrd}$
  \item $\sL$ is locally free outside a locus of codimension 3.
  \end{itemize}

  There is a sub-bundle
  \begin{equation*}
    \sV \incl \pls\sL
  \end{equation*}
  which, over each point in $\sgrd$, is equal to the aspect of the
  $\grd$ on the component containing the marked point.
\end{defn}

\begin{rmk}
  \label{rmk:pic-Mg1}
  By a theorem of Harer \cite{Harer}, for $g\ge3$,
  \begin{equation*}
    \Pic\Mtl\tensor\Q =
    \Q\lam\ds\Q\dn\ds\Q\de_1\ds\dotsb\ds\Q\de_{g-1}\ds\Q\psi
  \end{equation*}
  where $\lam$ and $\psi$ are the first Chern classes of the Hodge and
  tautological bundles respectively, $\dn$ is
  the divisor of irreducible nodal curves, and $\di$ is the divisor of
  unions of curves of genus $i$ and $g-i$, where the marked point lies
  on the component of genus $i$. 
\end{rmk}

\begin{defn}
  Since $\sL$ is locally free in codimension 3, the cycle classes
  \begin{equation*}
    c_1(\sL)^2 \cap [\scrd] \quad\text{ and }
    \quad c_1(\sL)\cdot c_1(\omega) \cap [\scrd]
  \end{equation*}
  on $\scrd$ are well-defined, where $\omega=\omega_{\scrd/\sgrd}$ is
  the relative dualizing sheaf. We define cycle classes
  $\ac,\bc,\cc\in A^1(\sgrd)$ as follows.
  \begin{align*}
    \ac &= \pls \bigl(c_1(\sL)^2 \cap [\scrd]\bigr) \\
    \bc &= \pls \bigl(c_1(\sL)\cdot c_1(\omega)\cap[\scrd]\bigr) \\
    \cc &= c_1(\sV)\cap[\sgrd]\text{.}
  \end{align*}
\end{defn}

We now state our main result.

\begin{thm}\label{main-thm}
  Let $g \ge 1$, $r \ge 0$, and $d \ge 1$ be integers for which
  \begin{equation*}
    \rho = g - (r+1)(g-d+r) = 0
  \end{equation*}
  and consider the map
  \begin{equation*}
    \eta\colon \sgrd ( \Mtl ) \to \Mtl
  \end{equation*}
  If
  \begin{equation*}
    \pf \colon \chow{1}{\sgrd (\Mtl)} \to \chow{1}{\Mtl}
  \end{equation*}
  is the proper push-forward on the corresponding Chow groups, then
  \begin{equation*}
    \begin{split}
      \frac{6(g-1)(g-2)}{dN}\pf\ac& = 
      6(gd - 2g^2 + 8d - 8g + 4) \lam \\
      & \quad + (2g^2 - gd + 3g - 4d - 2) \dn \\
      & \quad + 6\sum_{i=1}^{g-1} (g-i)(gd + 2ig - 2id - 2d) \di \\
      & \quad - 6d(g-2) \psi \text{,}
    \end{split}    
  \end{equation*}
  \begin{equation*}
    \begin{split}
      \frac{2(g-1)}{dN}\pf\bc
       = 12\lam - \dn 
        + 4 \sum_{i=1}^{g-1} (g-i)(g-i-1) \di  -  2(g-1)\psi\text{,}
    \end{split}
  \end{equation*}

  \begin{equation*}
    \begin{split}
      \frac{2(g-1)(g-2)}{N}\pf\cc
      & = \bigl[-(g+3)\mess + 5r(r+2)\bigr]\lam
      - d(r+1)(g-2)\psi \\
      & \quad + \frac{1}{6}\bigl[(g+1)\mess - 3r(r+2)\bigr]\dn \\
      & \quad + \sum_{i=1}^{g-1} (g-i) \bigl[i\mess +
      (g-i-2)r(r+2)\bigr]\di \text{,}
    \end{split}
  \end{equation*}
  where
  \begin{equation*}
    N = \frac{1!\cdot 2!\cdot 3! \cdots r!\cdot g!}
    {(g-d+r)!(g-d+r+1)!\cdots(g-d+2r)!}
  \end{equation*}
  and
  \begin{equation*}
    \mess = 3(g-1) + \frac{(r-1)(g+r+1)(3g-2d+r-3)}{g-d+2r+1} \text{.}
  \end{equation*}

\end{thm}

The proof of Theorem \ref{main-thm} will occupy the remaining two sections.

\section{Special Families of Curves}
\label{sec:spec-famil-curv}

Our strategy for proving Theorem \ref{main-thm} will be to pull back
to various families of stable curves over which the space of linear
series is easier to analyze. In Section \ref{sec:definitions-families}
we define three families of pointed curves, and in Section
\ref{sec:comp-spec-famil} we compute $\pf$ for these special families.
In Section \ref{sec:pull-back-maps} we compute the pull-backs of
the standard divisor classes on $\Mtl$ to the base spaces of each of
our families. Assembling the results of these three sections, we
compute $\pf$ over the whole moduli space.

\subsection{Definitions of Families}
\label{sec:definitions-families}

\begin{defn}
  \label{def:mogb}
  Let $i\colon \Mogb \incl \Mtl$ be the family of marked stable curves defined
  by sending a $g$-pointed stable curve
  \begin{equation*}
    (C,p_1,\dotsc,p_g)
  \end{equation*}
  of genus 0 to the stable curve
  \begin{equation*}
    \bigl( C \union \Union_{i=1}^g E_i , p_0 \bigr)
  \end{equation*}
  of genus $g$, where $E_i$ are fixed non-isomorphic elliptic curves,
  attached to $C$ at the points $p_i$, and $p_0\in E_1$ is fixed as
  well.
\end{defn}

\begin{figure*}[h]
  \unitlength 1mm
\begin{picture}(62.57,24.11)(0,0)
\linethickness{0.3mm}
\qbezier(3.65,5.84)(12.78,9.41)(19.28,9.49)
\qbezier(19.28,9.49)(25.78,9.56)(30.87,7.95)
\qbezier(30.87,7.95)(37.34,6.77)(48.25,8.02)
\qbezier(48.25,8.02)(59.18,9.27)(59.16,9.45)
\linethickness{0.3mm}
\qbezier(15.53,21.03)(12.66,16.49)(13.08,12.68)
\qbezier(13.08,12.68)(13.49,8.88)(14.67,4.79)
\linethickness{0.3mm}
\qbezier(25.03,20.88)(22.16,16.34)(22.58,12.53)
\qbezier(22.58,12.53)(22.99,8.73)(24.17,4.64)
\linethickness{0.3mm}
\qbezier(34.54,20.42)(31.67,15.88)(32.08,12.08)
\qbezier(32.08,12.08)(32.49,8.28)(33.67,4.19)
\linethickness{0.3mm}
\qbezier(53.98,20.58)(51.11,16.04)(51.53,12.23)
\qbezier(51.53,12.23)(51.94,8.43)(53.12,4.34)
\put(62.57,9.15){\makebox(0,0)[cc]{$C$}}

\put(15.23,24.05){\makebox(0,0)[cc]{$E_1$}}

\put(25.9,23.65){\makebox(0,0)[cc]{$E_2$}}

\put(35,23.4){\makebox(0,0)[cc]{$E_3$}}

\put(54.49,24.11){\makebox(0,0)[cc]{$E_g$}}

\put(45.06,23.89){\makebox(0,0)[cc]{$\dotsb$}}

\put(42.2,14.43){\makebox(0,0)[cc]{$\dotsb$}}

\put(14.46,18.24){\makebox(0,0)[cc]{$\bullet$}}

\put(10.68,18.45){\makebox(0,0)[cc]{$p_0$}}

\end{picture}
  \caption{$i(C,p_1,\dotsc,p_g)$}
\end{figure*}
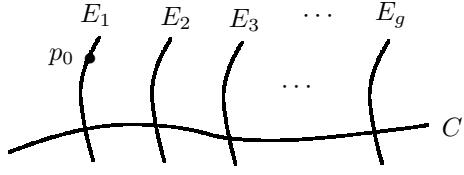

\begin{defn}
  \label{def:M21}
  Let $j\colon \Mttl \incl \Mtl$ be the family of curves defined by
  sending a marked curve $(C,p)$ to the marked stable curve
  \begin{equation*}
    (C \union C', p_0)
  \end{equation*}
  where $(C',p',p_0)$ is a fixed
  Brill-Noether-general
  curve in $\mtl_{g-2,2}$, attached nodally to $(C,p)$ at $p'$.
\end{defn}

\begin{figure*}[h]
%
%
%
%
%
%
\unitlength 1mm
\begin{picture}(58.24,13.11)(0,0)
\linethickness{0.3mm}
\qbezier(8.51,3.92)(17.95,3.67)(22.23,5.13)
\qbezier(22.23,5.13)(26.5,6.6)(30.61,12.83)
\linethickness{0.3mm}
\qbezier(19.41,13.11)(26.72,5.54)(33.67,3.95)
\qbezier(33.67,3.95)(40.62,2.35)(53.65,3.58)
\put(3.92,3.65){\makebox(0,0)[cc]{$C$}}

\put(58.24,3.38){\makebox(0,0)[cc]{$C'$}}

\put(45,3.11){\makebox(0,0)[cc]{$\bullet$}}

\put(45,6.35){\makebox(0,0)[cc]{$p_0$}}

\end{picture}
  \caption{$j(C,p)$}
\end{figure*}

\begin{defn}
  \label{def:marked-point}
  Fix Brill-Noether general curves
  \begin{align*}
    (C_1, p_1) & \in \m_{h,1} \\
    (C_2, p_2) & \in \m_{g-h,1} 
  \end{align*}
  and let $C = C_1 \union C_2$ be their nodal union along the $p_i$. Let
  $k_h\colon C_1 \incl \Mtl$ be the map sending $p\in C_1$ to the marked
  curve $(C,p)$.
\end{defn}

\subsection{Computations on the Special Families}
\label{sec:comp-spec-famil}

\begin{lemma}\label{lem:eta-mogb}
  For the family
  \begin{equation*}
    i\colon \Mogb \incl \Mtl
  \end{equation*}
  we have
  \begin{equation*}
    \pf\ac = \pf \bc = \pf\cc = 0
  \end{equation*}
\end{lemma}

\begin{lemma}\label{lem:eta-M21}
  For the family
  \begin{equation*}
    j \colon \Mttl \incl \Mtl
  \end{equation*}
  we have
  \begin{align*}
    \pf\ac
    &= \frac{2dN(d - 2g + 2)}{3(g-1)}(3\psi - \lam - \de_1)
    + \frac{dN}{g-1}(\lam + \de_1 - 4\psi) \\
    \pf\bc
    &= \frac{dN}{g-1}(\lam + \de_1 - 4\psi) \\
    \pf\cc
    &= \frac{-N\mess}{3(g-1)}(3\psi - \lam - \de_1) \text{,}
  \end{align*}
  where $N$ and $\mess$ are defined in the statement of Theorem \ref{main-thm}.
\end{lemma}

\begin{lemma}\label{lem:eta-marked-point}
  For the family
  \begin{equation*}
    k_h \colon C_1 \incl \Mtl
  \end{equation*}
  we have
  \begin{align*}
    \deg\pf\ac & = -d^2 N \\
    \deg\pf\bc & = -\bigl[2(g-h) - 1\bigr] dN \\
    \deg\pf\cc & = - \bigl[rh + \frac{1}{2}r(r+1)\bigr] N
  \end{align*}
\end{lemma}

\subsection{Pull-Back Maps on Divisors}
\label{sec:pull-back-maps}

\begin{lemma}
  \label{lem:pull-back-mogb}
  Let $\epsi$ be the class of the closure of the locus on $\Mogb$ of
  stable curves with two components, the component containing the
  first marked point having $i$ marked points.
  \begin{enumerate}
  \item The classes $\epsi$ are independent in $\coh{2}{\Mogb;\Q}$.
  \item
    For the family
    \begin{equation*}
      i\colon \Mogb \incl \Mtl
    \end{equation*}
    we have the following pull-back map on divisor classes.
    \begin{align*}
      i^* \lam & = i^*\psi = i^* \dn  = 0 \\
      i^* \di & = \epsi \qquad \textrm{for } i = 2,3,\ldots,g-2 \\
      i^*\de_1 & = -\sum_{i=2}^{g-2} \frac{(g-i)(g-i-1)}{(g-1)(g-2)}\epsi \\
      i^*\de_{g-1} & = -\sum_{i=2}^{g-2} \frac{(g-i)(i-1)}{g-2}\epsi \\
    \end{align*}
   \end{enumerate}
\end{lemma}

\begin{lemma}
  \label{lem:pull-back-M21}
  For the family
  \begin{equation*}
    j \colon \Mttl \incl \Mtl
  \end{equation*}
  we have the following pull-back map on divisor classes.
  \begin{align*}
    j^* \lam & = \lam & j^*\psi &= 0 \\
    j^* \dn & = \dn   & j^*\di &= 0 \quad i=1,2,\ldots,g-3 \\
    j^* \de_{g-2} & = -\psi & j^* \de_{g-1} & = \de_1
  \end{align*}
\end{lemma}

\begin{lemma}
  \label{pull-back-marked-point}
  For the family
  \begin{equation*}
    k_h\colon C_1 \incl \Mtl
  \end{equation*}
  we have the following pull-back map on divisor classes.
  \begin{align*}
    \deg k_h^*\lam & = 0 & \deg k_h^*\psi &= 2h -1 \\
    \deg k_h^*\de_h &= -1 & \deg k_h^*\de_{g-h} &= 1 \\
    \deg k_h^*\di &= 0 \quad i \ne h , g-h & & \\
  \end{align*}
\end{lemma}

\begin{proof}[Proof of Theorem \ref{main-thm}]
  Theorem \ref{main-thm} is now a consequence of the above lemmas.
  The main point is that the pull-backs of the classes
  $\pf\ac$, $\pf\bc$, and $\pf\cc$ to our special families coincide with the
  classes computed in Section~\ref{sec:comp-spec-famil}. For example,
  to see this for $j^*\pf\cc$, form the fiber the fiber square
  \begin{equation*}
    \begin{CD}
      j^*\sgrd @>{j'}>> \sgrd \\
      @V{\eta'}VV @VV{\eta}V \\
      \Mttl @>>j> \Mtl\text{.}
    \end{CD}
  \end{equation*}
  Notice that although $j$ is a regular embedding, $j'$ need not be.
  Nonetheless, according to Fulton \cite[Chapter 6]{Fulton.book}, there
  is a refined Gysin homomorphism
  \begin{equation*}
    j^! \colon A_k(\sgrd) \to A_{k-l}(j^*\sgrd)\text{,}
  \end{equation*}
  where $l$ is the codimension of $j$, which commutes with push-forward:
  \begin{equation*}
    \eta'_* j^! = j^*\pf\text{.}
  \end{equation*}
  We need to check that
  \begin{equation*}
    j^! c_1(\sV)\cap [\sgrd]
    = c_1({j'}^*\sV) \cap [j^*\sgrd]\text{.}
  \end{equation*}
  Since  
  \begin{equation*}
    j^! c_1(\sV)\cap [\sgrd]
    = c_1({j'}^*\sV) \cap j^![\sgrd]
  \end{equation*}
  \cite[Proposition~6.3]{Fulton.book}, it is enough to check that
  \begin{equation*}
    j^![\sgrd] = [j^*\sgrd]\text{.}
  \end{equation*}
  Generalizing the dimension upper bound in
  \cite[Corollary~5.9]{Osserman2} to the multi-component case 
  \cite{limit.linear.series}, we obtain
  \begin{equation*}
    \dim j^*\sgrd = \dim \Mttl \text{.}
  \end{equation*}
  This implies that the codimension of
  $j'$ is equal to that of $j$, so the normal cone of $j'$ is equal to
  the pull-back of the normal bundle of $j$, and the result follows.

  Now, for  example, to compute $\pf\cc$, write
  \begin{equation*}
    \pf\cc = a\lam - \sum_{i=0}^{g-1} b_i\di + c\psi
  \end{equation*}
  Our goal is to solve for $a,b_0,b_1,\dotsc,b_{g-1},c$. Using Lemmas
  \ref{lem:eta-marked-point} and \ref{pull-back-marked-point}, we may
  solve for $c$ and write $b_{g-i}$ in terms of $b_i$. From Lemmas
    \ref{lem:eta-mogb} and \ref{lem:pull-back-mogb}, we may further
    solve for $b_1,b_2,\dotsc,b_{g-2}$ in terms of $b_{g-1}$. It
    remains to determine $a$, $b_0$, and $b_{g-1}$. This is done by
    pulling back to $\Mttl$, which has Picard number 3, and using
    Lemmas \ref{lem:eta-M21} and \ref{lem:pull-back-M21}. The other
    push-forwards are computed similarly.
\end{proof}

\section{Proofs of Lemmas}
\label{sec:proofs-lemmas}

In this section, we give proofs of the lemmas stated in Sections
\ref{sec:comp-spec-famil} and \ref{sec:pull-back-maps}.

\begin{proof}[Proof of Lemma \ref{lem:eta-mogb}]
  \label{pf:eta-mogb}
  If $\cb_{0,g} \to \Mogb$ is the universal stable curve, then $i^*\Ctl$ is
  formed by attaching $\Mogb \times E_i$ to $\cb_{0,g}$ along the marked
  sections $\sigma_i \colon \Mogb \to \cb_{0,g}$. We have the
  following fiber square.
  \begin{equation*}
    \begin{CD}
      i^*\scrd @>>> i^*\Ctl \\
      @VVV @VVV \\
      i^*\sgrd @>>> \Mogb
    \end{CD}
  \end{equation*}

  By the Pl\"ucker formula for $\Proj^1$, given $[C]\in \Mogb$, a limit
  linear series on $i(C)$ must have the aspect
  \begin{equation*}
    (d-r-1)p_i + \ls{(r+1)p_i}
  \end{equation*}
  on each $E_i$. The line bundle $\sL \to i^*\scrd$ is, therefore,
  the pull-back from $i^*\Ctl$ of the bundle which is
  given by
  \begin{equation*}
    \pi_2^* \oh{E_1}{dp}
  \end{equation*}
  on $\Mogb \times E_1$ and is trivial on all other components. Thus
  $\ac=\bc=0$. The vector bundle $\sV \subset \pls \sL$ is trivial with
  fiber isomorphic to
  \begin{equation*}
    \coh{0}{\oh{E_1}{(r+1)p}}  %
    \subset
    \coh{0}{\oh{E_1}{dp}}
  \end{equation*}
  so $\cc=0$ as well.
\end{proof}

Before proving Lemma \ref{lem:eta-M21} we state an elementary result
in Schubert calculus.

\begin{lemma}\cite[p. 266]{GriffithsHarrisBN}
  \label{lem:schubert}
  For integers $r$ and $d$ with $0\le r\le d$, let
  \begin{equation*}
    X = \Grass(r,\Proj^d)
  \end{equation*}
  be the Grassmannian of $r$-planes in $\Proj^d$. For integers
  \begin{equation*}
    0 \le b_0 \le b_1 \le \dotsb \le b_r \le d-r \text{,}
  \end{equation*}
  let $\sch{b}=\sch{b_r,\dotsc,b_0}$ be the corresponding Schubert cycle of
  codimension $\sum b_i$. Let $\zeta = \sch{1,1,\dotsc,1,0}$ be the
  special Schubert cycle of codimension $r$.
  If $k$ is an integer for which
  \begin{equation*}
    rk + \sum_{i=0}^r b_i  = \dim X = (r+1)(d-r) \text{,}
  \end{equation*}
  then
  \begin{equation*}
    \int_X \zeta^k \cdot \sch{b}
    = \frac{k!}{\prod_{i=0}^r (k-d+r + a_i)! }
    \prod_{0\le i<j\le r}(a_j - a_i) \text{,}
  \end{equation*}
  where $a_i = b_i + i$.
\end{lemma}

\begin{proof}[Proof of Lemma \ref{lem:eta-M21}]
  \label{pf:eta-M21}
  Since $\mb_{2,1}$ is a smooth Deligne-Mumford stack, it is enough,
  by the moving lemma, to
  prove Lemma \ref{lem:eta-M21} for a family over a
  complete curve
  \begin{equation*}
    B \incl \Mttl
  \end{equation*}
  which intersects the boundary and Weierstrass divisors
  transversally.  If $\pi\colon\sC \to B$ is the universal stable
  genus-2 curve and $\sg\colon B\to \sC$ is the marked section, then
  $j^*\widetilde\sC_{g,1}$ is formed by attaching $\sC$ to $B \times
  C'$ along the marked section $\Sg=\sg(B) \subset \sC$.

  We begin by assuming that $B$ is disjoint from the closure of
  the Weierstrass locus $W$. In this case we claim that
  \begin{equation*}
    j^* \sgrd \to  B
  \end{equation*}
  is a trivial $N$-sheeted cover of the form $B\times X$, where is $X$
  a zero-dimensional scheme of length $N$.  Indeed, for any curve
  $(C,p)$ in $\Mttl \setminus W$ there are two (limit linear) $\grd$s
  on $C$ with maximum ramification at $p$; the vanishing sequences are
  \begin{align*}
    a_1 &= (d-r-2,d-r-1, \dotsc, d-4,d-3,d) \text{,} \\
    a_2 &= (d-r-2,d-r-1, \dotsc, d-4, d-2, d-1) \text{.}
  \end{align*}
  If $C$ is smooth, the two linear series are
  \begin{equation*}
    (d-r-2)p + \ls{(r+2)p}
  \end{equation*}
  and
  \begin{equation*}
    (d-r-2)p + \ls{rp + K_C} \text{.}
  \end{equation*}
  There are analogous series on nodal curves outside the closure of
  the Weierstrass locus. In the case of irreducible nodal curves, the
  sheaves are locally free.

  For each of the two $\grd$s on $C$ with maximum ramification at $p$,
  there are finitely many $\grd$s on $C'$ with compatible
  ramification. Specifically, there are
  \begin{equation*}
    \frac{(2g-2-d)N}{2(g-1)}
  \end{equation*}
  of type $a_1$ and
  \begin{equation*}
    \frac{dN}{2(g-1)} 
  \end{equation*}
  of type $a_2$, for total of $N$ limit linear series counted with
  multiplicity. Since $C'$ is
  fixed, the cover $j^*\sgrd \to B$ is a trivial $N$-sheeted cover.

  Consider a reduced sheet $B_1 \caniso B$ of type $a_1$. (We assume for
  simplicity that the sheet is reduced---the computation is the same
  in the general case.)  Then the universal
  line bundle $\sL$ on $j^*\sC^r_d$ is given as
  \begin{equation*}
    \sL \isom
    \begin{cases}
      \mathcal{O}_{\sC} & \text{on $\sC$} \\
      \pi_2^*L_1  & \text{on $B_1 \times C'$}
    \end{cases}
  \end{equation*}
  for some line bundle $L_1$ on $C'$ of degree $d$. It follows that
  $\ac=\bc=\cc=0$ on $B_1$.

  Next consider a sheet $B_2\caniso B$ of type $a_2$. Over $B_2 \times
  C'$ the universal line bundle $\sL$ is isomorphic to $\pi_2^* L_2$
  for some $L_2$ of degree $d$ on $C'$. It remains to determine $\sL$
  over $\sC$.
  Now $\omega_C(-2p)$ gives the correct line bundle for all $[C]\in
  B_2$; however, it has the wrong degrees on the components of the
  singular fibers. As our first approximation to $\sL$ on $\sC$ we take
  \begin{equation*}
    \omega_{\sC / B_2}( -2 \Sg)
  \end{equation*}
  Let $\De\subset\sC$ be the pull-back of
  the divisor on $\sC$ of curves of the form $C_1 \union C_2$, where
  the $C_i$ have genus one, and the marked points lie on different
  components. Then
  \begin{equation*}
    \omega_{\sC / B_2}( -2 \Sg + \De)
  \end{equation*}
  has the correct degree on the irreducible components on each
  fiber. It remains only to normalize our line bundle by pull-backs
  from the base $B_2$. In this case, $\sL|_\sC$ is required to be trivial
  along $\Sg$ since $\sL$ is a pull-back from $C'$ on the other
  component. If $\sg\colon B_2 \to \sC$  is the marked section, we
  let
  \begin{equation*}
    \Psi = \sg^* \omega_{\sC / B_2}
  \end{equation*}
  be the tautological line bundle on $B_2$. Then
  \begin{align*}
    \sg^* \oh{\sC}{\De} & \isom \mathcal{O}_{B_2} \\
    \sg^*\oh{\sC}{\Sg} & \isom \dual\Psi
  \end{align*}
  It follows that on $\sC$,
  \begin{equation*}
    \sL \isom
    \begin{cases}
      \omega_{\sC / B_2}( -2 \Sg + \De) \tensor \pi^*\Psi^{\tensor -3}
      &\text{on $\sC$} \\
      \pi_2^* L_2 &\text{on $B_2\times C'$.}
    \end{cases}
  \end{equation*}
  Thus, if we let
  \begin{align*}
    \omega &= c_1(\omega_{\sC / B_2}) \\
    \sg &= c_1 (\oh{\sC}{\Sg}) \\
    \de &= c_1 (\oh{\sC}{\De}) 
  \end{align*}
  on $\sC$ and let
  \begin{equation*}
    \psi = c_1(\Psi)
  \end{equation*}
  on $B_2$, then
  \begin{equation*}
    c_1(\sL) =
    \begin{cases}
      \omega - 2\sg + \de - 3\pi^*\psi  & \text{on $\sC$} \\
      d\pi_2^* p  & \text{on $B_2 \times C'$}
    \end{cases}
  \end{equation*}
  For the relative dualizing sheaf $\omega_{j^*\widetilde\sC_{g,1} /
    B_2}$, we have
  \begin{equation*}
    c_1(\omega_{j^*\sC / B_2}) =
    \begin{cases}
      \omega + \sg & \text{on $\sC$} \\
      (2(g-2) - 1)\pi_2^* p & \text{on $B_2 \times C'$.}
    \end{cases}
  \end{equation*}
  To compute the products of these classes on $j^*\sC^r_d$, recall the
  following formulas on $\pi\colon \overline{\sC}_{2,1} \to \mb_{2,1}$
  \begin{align*}
    \pls \omega & =2 &
    \pls \delta & = 0 &
    \pls \sigma & = 1 \\
    \pls \omega^2 & = 12\lam - \dn - \de_1 &
    \pls\sg^2 & = -\psi &
    \pls\de^2 &= -\de_1 \\
    \pls (\de . \sg) & = 0 &
    \pls (\omega . \de) & = \de_1 &
    \pls (\sg . \omega) & = \psi \text{.}
  \end{align*}
  Then we compute
  \begin{align*}
    \ac &= \pls\bigl[c_1(\sL)^2\bigr] = 12\lam - \dn - 8\psi \\
    \bc &= \pls\bigl[c_1(\sL) \cdot c_1(\omega)\bigr] = 12\lam - \dn - 8\psi
  \end{align*}
  on $B_2$.
  Since the marked point lies on
  $C'$, $\sV$ is trivial on $B_2$, so $\cc=0$ on $B_2$.

  \begin{figure*}
    \unitlength 1mm
\begin{picture}(68.38,56.35)(0,0)
\linethickness{0.3mm}
\qbezier(57.57,53.51)(55.24,47.78)(55.71,43.62)
\qbezier(55.71,43.62)(56.18,39.45)(58.11,35.14)
\qbezier(58.11,35.14)(59.89,30.84)(59.56,28.14)
\qbezier(59.56,28.14)(59.24,25.43)(57.36,22.3)
\linethickness{0.3mm}
\multiput(17.97,33.38)(0.12,0.45){43}{\line(0,1){0.45}}
\linethickness{0.3mm}
\multiput(17.43,42.84)(0.12,-0.39){54}{\line(0,-1){0.39}}
\linethickness{0.3mm}
\qbezier(47.43,52.03)(43.13,32.56)(39.33,35.44)
\qbezier(39.33,35.44)(35.53,38.32)(40.14,39.32)
\qbezier(40.14,39.32)(45.18,37.85)(46.17,31.78)
\qbezier(46.17,31.78)(47.17,25.71)(46.89,20.54)
\linethickness{0.3mm}
\qbezier(2.03,2.97)(12.44,6.23)(15.96,6.84)
\qbezier(15.96,6.84)(19.49,7.44)(32.65,3.92)
\qbezier(32.65,3.92)(47.82,1.03)(56.64,3.41)
\qbezier(56.64,3.41)(65.46,5.79)(60.42,4.32)
\linethickness{0.3mm}
\qbezier(2.43,24.73)(12.84,27.99)(16.37,28.6)
\qbezier(16.37,28.6)(19.89,29.2)(33.06,25.68)
\qbezier(33.06,25.68)(48.23,22.79)(57.05,25.17)
\qbezier(57.05,25.17)(65.86,27.55)(60.82,26.08)
\linethickness{0.3mm}
\qbezier(51.86,52.46)(49.54,46.72)(50,42.56)
\qbezier(50,42.56)(50.47,38.4)(52.41,34.08)
\qbezier(52.41,34.08)(54.19,29.79)(53.86,27.08)
\qbezier(53.86,27.08)(53.53,24.37)(51.66,21.24)
\linethickness{0.3mm}
\qbezier(34.95,53.03)(32.62,47.29)(33.09,43.13)
\qbezier(33.09,43.13)(33.55,38.96)(35.49,34.65)
\qbezier(35.49,34.65)(37.27,30.35)(36.94,27.64)
\qbezier(36.94,27.64)(36.62,24.94)(34.74,21.81)
\linethickness{0.3mm}
\qbezier(28.7,52.78)(26.38,47.05)(26.84,42.89)
\qbezier(26.84,42.89)(27.31,38.72)(29.24,34.41)
\qbezier(29.24,34.41)(31.02,30.11)(30.7,27.4)
\qbezier(30.7,27.4)(30.37,24.7)(28.5,21.57)
\linethickness{0.3mm}
\qbezier(13.68,51.73)(11.35,45.99)(11.82,41.83)
\qbezier(11.82,41.83)(12.28,37.67)(14.22,33.35)
\qbezier(14.22,33.35)(16,29.06)(15.67,26.35)
\qbezier(15.67,26.35)(15.35,23.64)(13.47,20.51)
\linethickness{0.3mm}
\qbezier(6.76,51.76)(4.43,46.02)(4.9,41.86)
\qbezier(4.9,41.86)(5.36,37.69)(7.3,33.38)
\qbezier(7.3,33.38)(9.08,29.08)(8.75,26.37)
\qbezier(8.75,26.37)(8.43,23.67)(6.55,20.54)
\put(23.11,56.35){\makebox(0,0)[cc]{$\Delta$}}

\put(66.89,26.62){\makebox(0,0)[cc]{$\Sigma$}}

\put(63.65,56.35){\makebox(0,0)[cc]{$\sC$}}

\put(67.7,5){\makebox(0,0)[cc]{$B_2$}}

\linethickness{0.3mm}
\put(40.68,7.16){\line(0,1){12.43}}
\put(40.68,7.16){\vector(0,-1){0.12}}
\end{picture}
    \caption{The morphism $\sC\to B_2$.}
  \end{figure*}
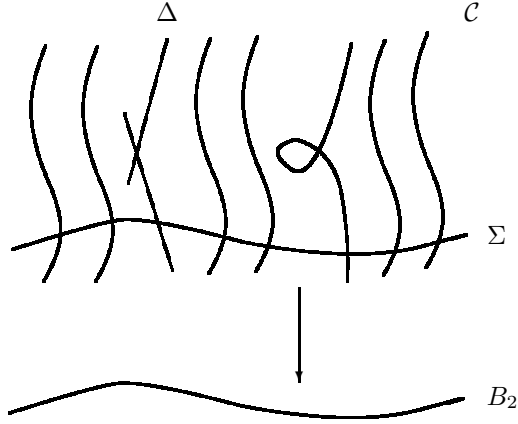

  Finally we consider the case where $B$ (transversally) intersects
  the Weierstrass locus. In this case
  \begin{equation*}
    \eta\colon j^*\sgrd \to B
  \end{equation*}
  is the union of a trivial $N$-sheeted cover of $B$ and a 
  1-dimensional scheme lying over each point of the divisor $W$. It
  will suffice to compute $\ac$ and $\cc$ on
  \begin{equation*}
    \Grd ( j[C,p])
  \end{equation*}
  where $C$ is a smooth genus-2 curve, and $p\in C$ is a Weierstrass
  point. (Note that $\bc$ is automatically zero.)

  There is a single
  $\grd$ on $C$ with maximal ramification at $p$, namely
  \begin{equation*}
    (d-r-2)p + \ls{(r+2)p}\text{,}
  \end{equation*}
  which has vanishing sequence
  \begin{equation*}
    (d-r-2,d-r-1,\dotsc,d-4,d-2,d)\text{.}
  \end{equation*}
  We claim that we only need to consider components of $\Grd(C\union
  C')$ with this aspect on $C$. Indeed, any $\grd$ on $C$ with
  ramification 1 less at $p$ is still a subseries of $\ls{dp}$. There
  will be finitely many corresponding aspects on $C'$ so that, as
  before, $\ac=\bc=\cc=0$ on these components.

  It remains to consider the components of $\Grd(C\union C')$ where
  the aspect on $C'$ has ramification $(0,1,2,2,\dotsc,2)$ or more at
  $p'$. We are reduced to studying the one-dimensional scheme
  \begin{equation*}
    S=\Grd (C'; p',(0,1,2,2,\dotsc,2)) \text{.}
  \end{equation*}
  To simplify computations we specialize $C'$ to a curve which is the
  union of $\Proj^1$ with $g-2$ elliptic curves $E_1,\dotsc,E_{g-2}$
  attached at general points $p_1,\dotsc,p_{g-2}$,
  and where the marked point $p_0$ lies on $E_1$, and the point of attachment
  $p'$ lies on the $\Proj^1$.
  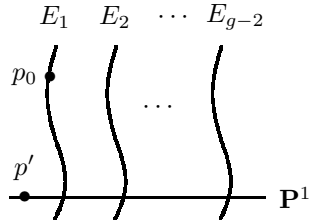
\begin{figure*}[h]
    \unitlength 1mm
\begin{picture}(42,28)(0,0)
\linethickness{0.3mm}
\put(4,4){\line(1,0){34}}
\linethickness{0.3mm}
\qbezier(10,24)(8.58,19.8)(8.86,16.76)
\qbezier(8.86,16.76)(9.15,13.71)(10.33,10.55)
\qbezier(10.33,10.55)(11.41,7.41)(11.21,5.43)
\qbezier(11.21,5.43)(11.02,3.45)(9.87,1.16)
\linethickness{0.3mm}
\qbezier(18,24)(16.58,19.8)(16.86,16.76)
\qbezier(16.86,16.76)(17.15,13.71)(18.33,10.55)
\qbezier(18.33,10.55)(19.41,7.41)(19.21,5.43)
\qbezier(19.21,5.43)(19.02,3.45)(17.87,1.16)
\linethickness{0.3mm}
\qbezier(32,24)(30.58,19.8)(30.86,16.76)
\qbezier(30.86,16.76)(31.15,13.71)(32.33,10.55)
\qbezier(32.33,10.55)(33.41,7.41)(33.21,5.43)
\qbezier(33.21,5.43)(33.02,3.45)(31.87,1.16)
\put(10,28){\makebox(0,0)[cc]{$E_1$}}

\put(18,28){\makebox(0,0)[cc]{$E_2$}}

\put(34,28){\makebox(0,0)[cc]{$E_{g-2}$}}

\put(24,16){\makebox(0,0)[cc]{$\dotsb$}}

\put(6,4){\makebox(0,0)[cc]{$\bullet$}}

\put(6,8){\makebox(0,0)[cc]{$p'$}}

\put(9.39,20){\makebox(0,0)[cc]{$\bullet$}}

\put(6,20){\makebox(0,0)[cc]{$p_0$}}

\put(42,4){\makebox(0,0)[cc]{$\Proj^1$}}

\put(26,28){\makebox(0,0)[cc]{$\dotsb$}}

\end{picture}
    \caption{The curve $C'$.}
  \end{figure*}
  There will be two types of components of $S$: those on which the
  aspects on the $E_i$ are maximally ramified at $p_i$, and those on
  which the aspect on one $E_i$ varies. Again, as in the proof of
  Lemma \ref{lem:eta-mogb}, we need only consider the latter case.

  Assume that for some $i$, the ramification at $p_i$ of the $\grd$ on
  $E_i$ is one less than maximal. There are two possibilities: either
  the series is of the form
  \begin{equation*}
    (d-r-1)p_i + \ls{rp_i+q} \qquad\text{for $q\in E_i$,}
  \end{equation*}
  which for $q\ne p_i$ imposes on the $\Proj^1$ the ramification condition
  \begin{equation*}
    (1,1,\dotsc,1) \text{,}
  \end{equation*}
  or the $\grd$ is a subseries of
  \begin{equation*}
    (d-r-2)p_i + \ls{(r+2)p_i}
  \end{equation*}
  containing
  \begin{equation*}
    (d-r)p_i + \ls{rp_i}\text{,}
  \end{equation*}
  which generically imposes on the $\Proj^1$ the ramification condition
  \begin{equation*}
    (0,1,1,\dotsc,1,2) \text{.}
  \end{equation*}
  In the first case the components are parameterized by $E_i$, and we
  compute that $\ac=-2$ on each such irreducible
  component, irrespective of whether $i=1$ or not. By
  Grothendieck-Riemann-Roch, $\cc=-1$ when $i=1$ and is zero otherwise. In
  the second case the $\grd$s are parameterized by a $\Proj^1$. Because
  the line bundle is constant, $\ac=0$. On each such $\Proj^1$,
  the vector bundle $\sV$
  may be viewed as the tautological bundle of rank $r+1$ on the
  Grassmannian of vector subspaces of a fixed vector space of
  dimension $r+2$ containing a subspace of dimension $r$. It follows
  that $\cc=-1$ on each $\Proj^1$.

  Let $X=\Grass(r,\Proj^d)$ be the Grassmannian of $r$-planes in $\Proj^d$. Let
  \begin{equation*}
    \zeta = \sch{1,1,\dotsc,1,0}
  \end{equation*}
  be the special Schubert cycle of codimension $r$.
  Collecting our calculations, we have that on $\Grd(C\union C')$,
  \begin{align*}
    \ac &= -2(g-2)\int_X
    \sch{2,2,\dotsc,2,1,0}\cdot\sch{1,1,\dotsc,1}\cdot\zeta^{g-3} \\ 
    &= -2(g-2)\int_X
    \sch{3,3,\dotsc,3,2,1}\cdot\zeta^{g-3} \text{,}
  \end{align*}
  and
  \begin{align*}
    \cc &= -\int_X \sch{2,2,\dotsc,2,1,0}\cdot
    (\sch{1,1,\dotsc,1} + \sch{2,1,1,\dotsc,1,0})\cdot
    \zeta^{g-3} \\
    &= -\int_X \sch{2,2,\dotsc,2,1,0}\cdot
    (\sch{1,0,0,\dotsc,0} \cdot \zeta)\cdot
    \zeta^{g-3} \\
    &= -\int_X
    (\sch{3,2,2,\dotsc,2,1,0}+\sch{2,2,\dotsc,2,0}+\sch{2,2,\dotsc,2,1,1})
    \cdot     \zeta^{g-2} \\ 
    &= -\int_X
    (\sch{3,2,2,\dotsc,2,1,0}+ \zeta^2)
    \cdot     \zeta^{g-2} \\ 
    &= -\int_X
    \sch{3,2,2,\dotsc,2,1,0} \cdot \zeta^{g-2}
    -\int_X \zeta^g \text{.}
  \end{align*}
  From Lemma \ref{lem:schubert} we compute,
  \begin{align*}
    \ac &= \frac{-2d(2g-2-d)N}{3(g-1)} \\
    \cc &= \frac{-\xi N}{3(g-1)} \text{.}
  \end{align*}
  Since the class of the Weierstrass locus in $\Mttl$ is $3\psi - \lam
  -\de_1$, the lemma follows.
\end{proof}

\begin{proof}[Proof of Lemma \ref{lem:eta-marked-point}]
  \label{pf:eta-marked-point}
  Because the curves $(C_i,p_i)$ are Brill-Noether general, $k_h^*
  \sgrd$ is a trivial $N$-sheeted cover of $C_1$ of the form
  $C_1\times X$, where $X$ is a zero-dimensional scheme of length $N$.
  Fix a sheet $G\isom C_1$ in $k_h^*\sgrd$; this choice corresponds to
  aspects
  \begin{equation*}
    V_i \subset \coh{0}{C_i, L_i}
  \end{equation*}
  where $L_i$ are degree-$d$ line bundles on $C_i$.
  If $(a_0,a_1,\dotsc,a_r)$ is the vanishing sequence of $V_1$ at
  $p_1$, then we know that
  \begin{align*}
    0 &= \rho(h,r,d) - \sum_{i=0}^r (a_i-i) \\
      &= (r+1)(d-r) -hr  - \sum_{i=0}^r a_i + \frac{1}{2}r(r+1)\text{,}
  \end{align*}
  so
  \begin{equation}
    \label{sum.a_i}
        \sum_{i=0}^r a_i = (r+1)d - \frac{1}{2}r(r+1) - hr \text{.}
  \end{equation}
  Let $\sC_1$ be the blow-up of $C_1 \times C_1$ at $(p_1,p_1)$, $E$
  the exceptional divisor, and $e$ its first Chern class. We may
  construct the universal curve $k_h^*\widetilde\sC_{g,1} \to C_1$ by
  attaching $C_1\times C_2$ to $\sC_1$ along $C_1 \times \{p_2\}$ and
  the proper transform of $C_1\times \{p_1\}$.  Over the sheet $G$,
  the universal line bundle $\sL$ on $k_h^*\sC^r_d$ is
  \begin{equation*}
    \pi_2^*L_1 \tensor \oh{\sC_1}{-dE}
    \tensor \pi_1^* \dual{L}_1(dp_1)
  \end{equation*}
  on $\sC_1$ and
  \begin{equation*}
    \pi_2^* L_2 (-dp_2) \tensor \pi_1^* \dual{L}_1
  \end{equation*}
  on $C_1 \times C_2$. Thus
  \begin{equation*}
    c_1(\sL) =
    \begin{cases}
      d \pi_2^* p - de & \text{on $\sC_1$} \\
      - d \pi_1^* p & \text{on $C_1 \times C_2$} \text{.}
    \end{cases}
  \end{equation*}
  The relative dualizing sheaf $\omega_{k_h^*\widetilde\sC_{g,1} /
    C_1}$ is isomorphic to
  \begin{equation*}
     \pi_2^* \omega_{C_1}\tensor \oh{\sC}{E} \tensor \pi_1^* \oh{C_1}{-p_1}
  \end{equation*}
  on $\sC_1$ and
  \begin{equation*}
    \pi_2^* \omega_{C_2}(p_2)
  \end{equation*}
  on $C_1 \times C_2$. We have
  \begin{equation*}
    c_1(\omega) =
    \begin{cases}
      - \pi_1^* p + (2h-2)\pi_2^* p + e & \text{on $\sC_1$} \\
      (2(g-h) - 1)\pi_2^* p & \text{on $C_1 \times C_2$} \text{.}
    \end{cases}
  \end{equation*}
  Thus, on $G$,
  \begin{align*}
    \deg\ac &= c_1(\sL)^2 = -d^2 \\
    \deg\bc &= c_1(\sL) \cdot c_1(\omega) = -d\bigl[2(g-j) - 1\bigr] \text{.}
  \end{align*}
  The formulas for $\pf\ac$ and $\pf\bc$ now follow.

  To calculate $\cc$ on $G$, notice that it suffices to
  compute $c_1(\sV')$, where
  \begin{equation*}
    \sV' = \sV \tensor L_1 (-d p_1)
  \end{equation*}
  is a sub-bundle of
  \begin{equation*}
    {\pi_1}_* (\pi_2^*L_1 (-dE)) \text{.}
  \end{equation*}
  We claim there is a bundle isomorphism
  \begin{equation*}
    \Ds_{j=0}^r\oh{C_1}{(a_j-d)p_1} \xrightarrow{\isom} \sV'
  \end{equation*}
  To describe the map, pick a basis $(\sg_0,\sg_1,\ldots,\sg_r)$ of
  $V_1\subset H^0(L_1)$ with $\sg_i$ vanishing to order $a_i$ at
  $p_1$. Given local sections $\tau_i$ of $\oh{C_1}{(a_i - d)p_1}$, let the
  image of $(\tau_0,\tau_1,\dotsc,\tau_r)$ be the section
  \begin{equation*}
    \sum_{i=0}^r \sg_i \tau_i
  \end{equation*}
  of $\sV'$. This is clearly an isomorphism away from $p_1$ and is
  checked to be an isomorphism over $p_1$ as well. 
  Using \eqref{sum.a_i}, we have that on $G$,
  \begin{equation*}
    \begin{split}
      \deg\cc = \deg\sV' &= \sum_{i=0}^r (a_i-d) \\
      &= -\frac{1}{2}r(r+1) - rh
    \end{split}
  \end{equation*}
  which finishes the proof of the lemma.
\end{proof}

\begin{proof}[Proof of Lemma \ref{lem:pull-back-mogb}]
  \label{pf:pull-back-mogb}
  To prove the independence of the $\epsi$, consider the curves
  \begin{equation*}
    B_j \incl \Mogb
  \end{equation*}
  for $j=2,3,\dotsc,g-3$ given by taking a fixed stable curve in
  $\eps_j$ and moving a marked point on the component with $g-j$
  marked points. Let $B_1\incl\Mogb$ be the curve given by moving the
  first marked point along a fixed smooth curve.%
  The intersection matrix
  \begin{equation*}
    (\epsi \cdot B_j)
  \end{equation*}
  is
  \begin{equation*}
    \begin{pmatrix}
      g-1 &0&0&0&\cdots&0&0&0 \\
      -1&1&0&0&\cdots&0&0&g-3 \\
      0&-1&1&0&\cdots&0&0&g-4 \\
      \vdots& & & &\ddots& & & \vdots \\
      0&0&0&0&\cdots&-1&1&3 \\
      0&0&0&0&\cdots&0&-1&2 
    \end{pmatrix}
  \end{equation*}
  where the rows correspond to the $B_j$ for $j=1,2,\dotsc,g-3$, and
  the columns to the $\epsi$ for $i=2,3,\dotsc,g-2$. Since this matrix
  is non-singular, the first part of the lemma follows.

  To derive the formula for the pull-back, we follow Harris and
  Morrison \cite[Section 6.F]{HarMor}. Let $B$ be a smooth projective
  curve, $\pi\colon\sC \to B$ a 1-parameter family of curves in
  $\Mogb$ transverse to the boundary strata. Then $\pi$ has smooth
  total space, and the fibers of $\pi$ have at most two irreducible
  components. Let
  \begin{equation*}
    \sg_i \colon B\to\sC
  \end{equation*}
  be the marked sections. Denote by $\Sg_i$ the image curve $\sg_i(B)$
  in $\sC$. Then on $B$,
  \begin{align*}
    \de_1 &= \Sg_1^2 \\
    \de_{g-1} &= \sum_{j=2}^g \Sg_j^2 \text{,}
  \end{align*}
  where we are using $D^2$ to denote $\pls(D^2)$ for a divisor $D$ on $\sC$.
  We now contract the component of each reducible fiber which meets the
  section $\Sg_1$. If $\Sgt_j$ is the image of $\Sg_j$ under this
  contraction, then we have
  \begin{align*}
    \Sg_1^2 &= \Sgt_1^2 - \sum_{i=2}^{g-2}\epsi \\
    \sum_{j=2}^g \Sg_j^2 &= \sum_{j=2}^g\Sgt^2_j - \sum_{i=2}^{g-2}(i-1)\epsi
  \end{align*}
  The $\Sgt_j$ are sections of a $\Proj^1$-bundle, so
  \begin{equation*}
    \begin{split}
      0 &= (\Sgt_j - \Sgt_k)^2 
        = \Sgt_j^2 + \Sgt_k^2 - 2\Sgt_j\cdot\Sgt_k
    \end{split}
  \end{equation*}
  Thus
  \begin{equation*}
    \begin{split}
      (g-2)\sum_{j=2}^g\Sgt_j^2
      &= \sum_{2\le j,k\le g} (\Sgt_j^2 + \Sgt_k^2) 
       = 2\sum_{2\le j,k\le g} \Sgt_j\cdot\Sgt_k \\
      &= 2\sum_{i=2}^{g-2}\binom{i-1}{2}\epsi \text{.}
    \end{split}
  \end{equation*}
  It follows that
  \begin{equation*}
    \begin{split}
      \de_{g-1}
      &=
      \sum_{i=2}^{g-2}\left[\frac{(i-1)(i-2)}{g-2}-(i-1)\right]\epsi
      = \sum_{i=2}^{g-2}\frac{(i-1)(i-g)}{g-2}\epsi
    \end{split}
  \end{equation*}
  Similarly, we can show that
  \begin{equation*}
    \de_1 + \de_{g-1}
    = \sum_{i=2}^{g-2}\frac{i(i-g)}{g-1}\epsi
  \end{equation*}
  so the formula for $\de_1$ follows as well.
\end{proof}

The proofs of Lemmas \ref{lem:pull-back-M21} and
\ref{pull-back-marked-point} are straightforward, so we omit them.

\section{Appendix}
\label{sec:appendix}

\begin{prop}
  Let $\Sigma$ be a set of 21 general points in $\Proj^2$ and let $S =
  \Bl_{\Sigma}\Proj^2$ be the blow-up of $\Proj^2$ at $\Sigma$. If $H$
  is the line class on $S$ and $E_1,\dotsc,E_{21}$ are the exceptional
  divisors, then the linear system
  \begin{equation*}
     \Bigl\lvert 13H - 2\sum_{j=1}^9E_j - 3\sum_{k=10}^{21}E_k \Bigr\rvert
  \end{equation*}
  on $S$ contains a smooth connected curve.
\end{prop}

\begin{proof}
  We begin by showing that it is enough to exhibit a single set of 21
  points over a finite field for which the above statement is true.

  Let $\hilb{k}$ be the Hilbert scheme of $k$ points in $\Proj^2$, and
  let
  \begin{equation*}
    \Sigma_k \subset \hilb{k}\times\Proj^2
  \end{equation*}
  be the universal
  subscheme. Let
  \begin{equation*}
    B\subset\hilb{9}\times\hilb{12}
  \end{equation*}
  be the irreducible
  open subset over which the composition
  \begin{equation*}
    \begin{CD}
      \Sigma = \pi_1^{-1}\Sigma_9 \union \pi_2^{-1}\Sigma_{12}
    @>>> \hilb{9}\times\hilb{12} \times \Proj^2 \\
    @. @VVV \\
     @. \hilb{9}\times\hilb{12} \\
    \end{CD}
  \end{equation*}
  is \'etale, where $\pi_1$ and $\pi_2$ are the obvious projections. Let
  \begin{equation*}
    \pi\colon S = \Bl_\Sigma\Proj^2_B \to B
  \end{equation*}
  be the smooth surface over B whose fibers are blow-ups of $\Proj^2$
  at 21 distinct points. If $E_9$ and $E_{12}$ are the exceptional
  divisors, let
  \begin{equation*}
    \sL = \oh{S}{13H - 2E_9 -3E_{12}}\text{.}
  \end{equation*}
  We may further restrict $B$ to an open over which $\pi_*\sL$ is
  locally free of rank at least 6.

  If
  \begin{equation*}
    \sC \subset \Proj\pi_*\sL \times_B S
  \end{equation*}
  is the universal section, then the projection
  \begin{equation*}
    \sC \to \Proj\pi_*\sL
  \end{equation*}
  is flat, so it suffices to find a single smooth fiber in order to
  conclude that the general fiber is smooth. To this end we use
  Macaulay 2 \cite{Macaulay2} and work over a finite field.

\begin{verbatim}

i1 : S = ZZ/137[x,y,z];

\end{verbatim}
  Following Shreyer and Tonoli \cite{Schreyer.Tonoli}, we realize our
  points in $\Proj^2$ as a determinental subscheme.
\begin{verbatim}

i2 : randomPlanePoints = (delta,R) -> (
       k:=ceiling((-3+sqrt(9.0+8*delta))/2);
       eps:=delta-binomial(k+1,2);
       if k-2*eps>=0
       then minors(k-eps,
         random(R^(k+1-eps),R^{k-2*eps:-1,eps:-2}))
       else minors(eps,
         random(R^{k+1-eps:0,2*eps-k:-1},R^{eps:-2})));

i3 : distinctPoints = (J) -> (
       singJ = minors(2, jacobian J) + J;
       codim singJ == 3);

\end{verbatim}
  Let $\Sigma_9$ and $\Sigma_{12}$ be our subsets of 9 and 12 points,
  respectively.
\begin{verbatim}

i4 : Sigma9 = randomPlanePoints(9,S);

o4 : Ideal of S

i5 : Sigma12 = randomPlanePoints(12,S);

o5 : Ideal of S

i6 : (distinctPoints Sigma9, distinctPoints Sigma12)

o6 = (true, true)

o6 : Sequence

\end{verbatim}
  Their union is $\Sigma$.
\begin{verbatim}

i7 : Sigma = intersect(Sigma9, Sigma12);

o7 : Ideal of S

i8 : degree Sigma

o8 = 21

\end{verbatim}
  Next we construct the 0-dimensional subscheme $\Gamma$ whose ideal
  consists of curves double through points of $\Sigma_9$ and triple
  through points of $\Sigma_{12}$.
\begin{verbatim}

i9 : Gamma = saturate intersect(Sigma9^2, Sigma12^3);

o9 : Ideal of S

\end{verbatim}
  Let us check that $\Gamma$ imposes the expected number of conditions
  ($9\cdot 3+12\cdot 6=99$) on curves of degree 13.
\begin{verbatim}

i10 : hilbertFunction (13, Gamma)

o10 = 99

\end{verbatim}
  Pick a random curve $C$ of degree 13 in the ideal of $\Gamma$.
\begin{verbatim}

i11 : C = ideal (gens Gamma 
                  * random(source gens Gamma, S^{-13}));

o11 : Ideal of S

\end{verbatim}
  We check that $C$ is irreducible.
\begin{verbatim}

i12 : # decompose C

o12 = 1

\end{verbatim}
  \newcommand{\Csing}{C_{\mathrm{sing}}}%
  To check smoothness, let $\Csing$ be the singular locus of $C$.
\begin{verbatim}

i13 : Csing = (ideal jacobian C) + C;

o13 : Ideal of S

i14 : codim Csing

o14 = 2

\end{verbatim}

  A double point will contribute 1 to the degree of $\Csing$ if it is
  transverse and more otherwise. Similarly, a triple point will
  contribute 4 to the degree of $\Csing$ if it is transverse and more
  otherwise. So for $C$ to be smooth in the blow-up, we must have that
  \begin{equation*}
    \deg\Csing = 9 + 4\cdot 12 = 57
  \end{equation*}
\begin{verbatim}

i15 : degree Csing

o15 = 57

\end{verbatim}
\end{proof}

\providecommand{\bysame}{\leavevmode\hbox to3em{\hrulefill}\thinspace}
\providecommand{\href}[2]{#2}

\end{document}